\documentclass[11pt]{amsart}

\usepackage{fullpage, amsmath, amssymb, amsthm, enumerate, url, bibentry, hyperref, color,xcolor,mathtools}
\usepackage{todonotes, float,makecell,multirow}
\usepackage{graphicx} 
\usepackage{subcaption} 
\usepackage{booktabs}   
\usepackage{multirow}   
\usepackage{caption}   
\usepackage{algorithm}
\usepackage{algpseudocode}
\usepackage[normalem]{ulem}

\usepackage{tikz}
\usepackage{bm} 
\usetikzlibrary{shapes,arrows,shadows,calc,fit,decorations.pathreplacing}

\pgfdeclarelayer{background}
\pgfdeclarelayer{midground} 
\pgfdeclarelayer{foreground}
\pgfsetlayers{background,midground,main,foreground}

\tikzstyle{sensor}=[draw, fill=blue!20, text width=5em, text centered, minimum height=2.5em,drop shadow]
\tikzstyle{tinysensor}=[draw, fill=blue!20, text width=2em, text centered, minimum height=1.5em]
\tikzstyle{tinysign}=[draw, fill=green!20, text width=12.5em, text centered , minimum height=1.5em]
\tikzstyle{ann} = [above, text width=5em, text centered]
\tikzstyle{wa} = [sensor, text width=12em, fill=red!20, minimum height=10em, rounded corners, drop shadow]
\tikzstyle{sc} = [sensor, text width=13em, fill=red!20, minimum height=10em, rounded corners, drop shadow]
\tikzstyle{decision} = [diamond, draw, fill=orange!20, text width=6em, text badly centered, node distance=3cm, inner sep=0pt, drop shadow]
\tikzstyle{bubble}=[circle, draw=black!50, minimum size=20pt, inner sep=0pt]
\tikzset{mybox/.style={draw, fill=white, text width=13em, align=center,minimum height=1.5em, font=\small, inner xsep=0pt}}
\tikzset{smallbox/.style={draw, fill=white, text width=12em, align=center,minimum height=1.5em, font=\small, inner xsep=0pt}}
\tikzset{block arrow/.style={single arrow,draw=black,fill=black!70,line width=1pt,minimum height=2.5em, minimum width=1em, single arrow head extend=1em, font=\color{white}, align=center}}

\numberwithin{theorem}{section}
\numberwithin{equation}{section}
\numberwithin{table}{section}

\def\sgn{{\rm sgn}}
\def\C{\mathbb C}
\def\F{\mathbb F}
\def\Z{{\mathbb Z}}

\def\Q{{\mathbb Q}}

\def\sgn{{{\operatorname{sgn}}}}

\def\freq{{{\operatorname{freq}}}}
\def\MCC{{{\operatorname{MCC}}}}
\def\ECQ6{{{\operatorname{ECQ6}}}}
\def\ECQ8{{{\operatorname{ECQ8}}}}

\title{Twist Class Redundancy Drives the Prediction of Traces of Frobenius of Elliptic Curves}

\author{Angelica Babei}
\address{
  Department of Mathematics,
  Howard University,
  204 Academic Service Building B
  Washington, D.C. 20059
  USA
}
\email{\href{mailto:angelica.babei@howard.edu}{angelica.babei@howard.edu}}
\urladdr{\url{https://angelicababei.com}}

\author{Ujjawal Shah}
\address{
  Department of Electrical Engineering and Computer Science,
  Howard University,
  Washington, D.C. 20059
  USA
}
\email{\href{ujjawal.shah@bison.howard.edu}{ujjawal.shah@bison.howard.edu}}

\author{Malick Kebe}
\address{
  Department of Mathematics,
  Howard University,
  204 Academic Service Building  B
  Washington, D.C. 20059
  USA
  }
\email{\href{malick.kebe@bison.howard.edu}{malick.kebe@bison.howard.edu}}

\begin{document}

\begin{abstract}
Recent interest in applying machine learning methods to predict invariants of mathematical objects has yielded models with surprisingly strong performance, including those predicting traces of Frobenius for elliptic curves. We demonstrate that the underlying datasets contain significant redundancy within quadratic twist classes, which alone is sufficient to produce highly accurate predictions. To ensure future models capture new arithmetic properties rather than potentially exploiting these dataset artifacts, we introduce a benchmark dataset consisting exclusively of unique twist class representatives.
\end{abstract}

\maketitle

\section{Introduction}\label{sec:intro}
The advent of large databases of arithmetic objects, most notably the $L$-functions and Modular Forms Database (LMFDB) \cite{LMFDB}, has naturally led number theorists to ask whether machine learning (ML) can aid in the discovery of new mathematical phenomena. A significant success in this direction was the discovery of ``murmurations" of elliptic curves, which are unexpected oscillatory patterns in the statistics of Frobenius traces, and which were first identified through ML-driven observation \cite{HLOP24}. These efforts sparked intense activity, resulting in more than ten papers within two years that demonstrate murmurations across a variety of motivic objects, including modular forms, Maass forms, and Dirichlet characters \cite{BBLLD, BLLD+, BKM, Cowan, C1, LOP25, Martin, SawinSutherland2025, Zubrilina}. Beyond murmurations, the use of ML has extended more broadly across number theory and arithmetic geometry. For example, recent studies have applied neural-network and tree-based models to predict ranks of elliptic curves from traces of Frobenius \cite{KV23}, to classify class groups \cite{AHLOS} and Galois groups \cite{KS25} of number fields, to predict the order of the Shafarevich-Tate group \cite{BBSHS24}, or to predict traces of Frobenius of elliptic curves \cite{BCCHLLNP}.

In \cite{BCCHLLNP}, the authors investigate whether one can predict traces of Frobenius $a_p(E)$ for elliptic curves of bounded conductor based on nearby traces $(a_q(E))_{q<100, q\ne p}$ for elliptic curves of bounded conductor $N(E) < 10^6$. The study used transformer- and neural-network based models, and attained high Matthews Correlation Coefficients (MCC) for the magnitude $|a_p(E)|$ and for local properties such as $a_p(E) \bmod 2$. While these results suggest that transformers can identify patterns in traces of Frobenius, the ``black-box" nature of these models makes it difficult to determine whether they are learning deep arithmetic properties or simply exploiting statistical redundancies within the datasets.

A source of such redundancy arises from families of quadratic twists. If \(E\) and \(E^\chi\) are quadratic twists, then for all but finitely many primes \(p\),
\[
a_p(E^\chi)=\chi(p)a_p(E),
\]
and hence
\[
|a_p(E^\chi)|=|a_p(E)|.
\]
A model that learns to identify the quadratic twist family of a curve can reduce the prediction problem to recovering signs from previously observed data. Motivated by this observation, we propose  an algorithm that uses the absolute values of a small collection of traces of Frobenius
\[
(|a_{p_1}(E)|,\dots,|a_{p_k}(E)|)
\]
as a proxy for identifying the quadratic twist family of \(E\). The algorithm then predicts \(a_p(E)\) by determining whether the training set contains a quadratic twist of \(E\). We use this algorithm as a baseline test for the dataset of elliptic curves with $N(E) < 10^6$. While the original transformer models reached a maximum MCC of at most $0.6$ when predicting exact values, our algorithm significantly outperforms previously trained models, achieving an MCC of $0.79$.

This suggests that much of the available predictive power in these datasets can be captured through twist identification. To move toward detecting more subtle, non-trivial arithmetic properties, we introduce a refined dataset \cite{ECQ7} that accounts for these redundancies.

\section{Preliminaries}
\label{sec: prelim}

\subsection{Elliptic curves} \label{subsec: ellcurves} 
We recall a few preliminaries on elliptic curves. This material is standard, for example see \cite{silvermanArithmetic}. An elliptic curve $E$ over $\mathbb{Q}$ is a smooth, projective curve  of genus 1 with a distinguished point at infinity, and can be described by a Weierstrass equation of the form
\begin{equation} 
\label{eq: ell}
E : y^2 = x^3 + Ax + B, \qquad A, B \in \mathbb{Q},
\end{equation}

$E$ is said to have \textit{good reduction} at a prime $p$ if the equation \ref{eq: ell} defining $E$ over $\F_p$ is nonsingular. $E$ has good reduction at all but finitely many primes. At such (finitely many) other primes, $E$ is said to have \textit{bad reduction}. The conductor $N(E)$ is a product of powers of these primes with bad reduction, with exponents being dictated by whether $E$ has multiplicative or additive reduction.

We can associate to $E$ an  $L$-function with Euler product
$$L(E, s) = \prod_p L_p(E, s).$$
where if $E$ has good reduction at $p$, the local factor, sometimes referred to as an Euler factor, is
$$L_p(E, s) = \frac{1}{1 - a_p(E)p^{-s} + p^{1-2s}}.$$

For a prime $\ell$, we have a 2-dimensional representation $\rho_{E, \ell}:\operatorname{Gal}(\overline{\mathbb{Q}}/\mathbb{Q}) \to \operatorname{GL}_2(\mathbb{Z}_\ell)$ of the absolute Galois group $\operatorname{Gal}(\overline{\mathbb{Q}}/\mathbb{Q})$.
At primes $p \neq \ell$ where $E$ has good reduction, the coefficient $a_p(E)$ is the \textit{trace of the (arithmetic) Frobenius} element $\operatorname{Frob}_p$ under this representation
$$a_p(E) = \operatorname{tr}(\rho_{E, \ell}(\operatorname{Frob}_p)),$$ and is independent of $\ell$. In fact, one has in such cases that $$a_p(E) = p+1-|E(\F_p)|,$$ where $|E(\F_p)|$ is the number of solutions of $E$ over the finite field $\F_p$.

 An isogeny between two elliptic curves is a non-constant rational map that is also a homomorphism of algebraic groups. Being isogenous is an equivalence relation, and the set of all elliptic curves over $\mathbb{Q}$ isogenous to $E$ makes up its \textit{isogeny class}. Two elliptic curves over $\mathbb{Q}$ are isogenous if and only if they share the same $L$-function, meaning they have the exact same sequence $(a_p)_{p \text{ prime}}$.

Given a quadratic Dirichlet character $\chi$, the \textit{quadratic twist} $E^\chi$ is an elliptic curve isomorphic to $E$ over the quadratic field corresponding to $\chi$. For all primes $p$ of good reduction that do not divide the conductor of $\chi$, the trace of Frobenius of the twist is given by
$$a_p(E^\chi) = \chi(p)a_p(E).$$
Therefore, two quadratic twists share the exact same trace magnitudes, $|a_p(E^\chi)| = |a_p(E)|$, at all but finitely many primes.

In order to differentiate between curves in distinct twist classes, we use a \textit{twist hash} as defined in Sec 4.3 of \cite{BSSVY} and implemented in Magma \cite{Magma, SutherlandMagma}. Specifically, given an elliptic curve $E$ over $\Q$, the twist hash of $E$ is calculated as
\[ h(E) \equiv \sum_{2^{12} < p < 2^{13}} c_p |a_p(E)| \pmod{P},
\]
where $P=2^{61}-1$ and $c_p = \left\lfloor \pi\, P^{e_p} \right\rfloor$ with $e_p := \#\{q : 2^{12} < q \le p,\ q \text{ prime}\},$ and $\pi=3.14 \dots$. Here, $c_p$ corresponds to digits of $\pi$ in base $P$. In datasets of size $ \ll 2^{61}$, hash collisions are highly unlikely.

\subsection{Clustering metrics} 
\label{subsec: cluster}

In section \ref{sec: twists}, we use a proxy for partitioning the curves into twist classes by the absolute values of the last $k$ primes of the input $(|a_{p_1}|, |a_{p_2}|, \dots, |a_{p_k}|)$. To measure how well this partition matches the true partition defined by the twist hash, we use standard external clustering metrics. Let
\[
\Pi^{\mathrm{true}}=\{C_1,\dots,C_r\}
\]
denote the true twist-hash partition and let
\[
\Pi^{(k)}=\{K_1,\dots,K_s\}
\]
denote the partition induced by the last \(k\) primes. Write
\[
n_{ij}=|C_i\cap K_j|,\qquad
a_i=\sum_{j=1}^s n_{ij},\qquad
b_j=\sum_{i=1}^r n_{ij},
\]
where \(n\) is the total number of curves.

Our primary metric is the \emph{Adjusted Rand Index} (ARI), which compares the two partitions at the level of pairs of observations \cite{HubertArabie1985}. In terms of the contingency table \((n_{ij})\), ARI is
\[
\mathrm{ARI}(\Pi^{\mathrm{true}},\Pi^{(k)})
=
\frac{
\displaystyle \sum_{i,j}\binom{n_{ij}}{2}
-
\frac{
\displaystyle \left(\sum_i \binom{a_i}{2}\right)\left(\sum_j \binom{b_j}{2}\right)
}{
\displaystyle \binom{n}{2}
}
}{
\displaystyle
\frac12\left[
\sum_i \binom{a_i}{2}
+
\sum_j \binom{b_j}{2}
\right]
-
\frac{
\displaystyle \left(\sum_i \binom{a_i}{2}\right)\left(\sum_j \binom{b_j}{2}\right)
}{
\displaystyle \binom{n}{2}
}
}.
\]
This is our main criterion for selecting \(k\), since it penalizes both accidental mergers of distinct twist-hash classes and oversplitting of repeated classes.

We also report \emph{homogeneity}, \emph{completeness}, and \emph{V-measure}, following Rosenberg and Hirschberg \cite{RosenbergHirschberg2007}. If \(C\) denotes the true class label and \(K\) the predicted cluster label, then
\[
H(C)=-\sum_{i=1}^r \frac{a_i}{n}\log\!\left(\frac{a_i}{n}\right),
\qquad
H(K)=-\sum_{j=1}^s \frac{b_j}{n}\log\!\left(\frac{b_j}{n}\right),
\]
and
\[
H(C\mid K)
=
-\sum_{j=1}^s\sum_{i=1}^r
\frac{n_{ij}}{n}\log\!\left(\frac{n_{ij}}{b_j}\right),
\qquad
H(K\mid C)
=
-\sum_{i=1}^r\sum_{j=1}^s
\frac{n_{ij}}{n}\log\!\left(\frac{n_{ij}}{a_i}\right),
\]
with the usual convention that zero terms are omitted. Homogeneity and completeness are then defined by
\[
h=
\begin{cases}
1, & \text{if } H(C)=0,\\[4pt]
1-\dfrac{H(C\mid K)}{H(C)}, & \text{otherwise},
\end{cases}
\qquad
c=
\begin{cases}
1, & \text{if } H(K)=0,\\[4pt]
1-\dfrac{H(K\mid C)}{H(K)}, & \text{otherwise},
\end{cases}
\]
and the V-measure is their harmonic mean,
\[
V=
\begin{cases}
0, & \text{if } h+c=0,\\[4pt]
\dfrac{2hc}{h+c}, & \text{otherwise}.
\end{cases}
\]

\subsection{Comparison metrics} 
\label{subsec: comparison}
We further introduce three metrics to assess predictive accuracy and to quantify agreement between models. To evaluate predictive accuracy over a test set of size $N$ with true values $y_i$ and predicted values $\hat{y}_i$, we utilize Mean Absolute Error (MAE) and Root Mean Squared Error (RMSE), where
$$ \text{MAE} = \frac{1}{N} \sum_{i=1}^{N} |y_i - \hat{y}_i|, \quad \text{RMSE} = \sqrt{\frac{1}{N} \sum_{i=1}^{N} (y_i - \hat{y}_i)^2} $$

Furthermore, to quantify the divergence between two distinct models (e.g., say Model A and Model B) yielding predictions $\hat{y}_i^{(A)}$ and $\hat{y}_i^{(B)}$ on the same test set, we define \textit{disagreement} as the proportion of mismatched predictions:
$$ \text{Disagreement} = \frac{1}{N} \sum_{i=1}^{N} \mathbf{1}\left(\hat{y}_i^{(A)} \neq \hat{y}_i^{(B)}\right) $$
where $\mathbf{1}$ denotes the indicator function yielding $1$ when the condition is true and $0$ otherwise.

\subsection{Performance metrics}
\label{s: mcc}

In addition to accuracy, we utilize the Mathews correlation coefficient (MCC), especially in classification experiments with unbalanced classes. For binary classification experiments, the MCC \cite{matthews1975comparison} is defined as 
$$\MCC = \frac{TP \times TN - FP \times FN}{\sqrt{(TP + FP)(TP + FN)(TN + FP)(TN + FN)}},$$
where $TP, FP, TN$ and $FN$ are the number of true positives, false positives, true negatives and false negatives, respectively. The MCC attains values in the range $[-1, 1]$, where $1$ indicates perfect prediction, $0$ indicates no correlation between predictions and true values, and $-1$ indicates perfectly wrong classification.

In \cite{gorodkin2004comparing}, the MCC is generalized for multiclass experiments by

$$\MCC = \frac{\sum_k \sum_l \sum_m C_{kk}C_{lm} - C_{kl}C_{mk}}{\sqrt{\sum_k 
 (\sum_l C_{kl}) (\sum_{k'| k' \ne k} \sum_{l'} C_{k'l'})} \sqrt{\sum_k 
 (\sum_l C_{lk}) (\sum_{k'| k' \ne k} \sum_{l'} C_{l'k'})} },$$
 where $C_{ij}$ is the $(i,j)$ entry in the confusion matrix $C$. 

Unlike accuracy, MCC remains informative under class imbalance, which is essential in our setting because the empirical distribution of $a_p(E)$ is highly non-uniform.

\subsection{Notation and conventions} \label{subsec: notation}

For a real number $x$, we define the sign function $\operatorname{sgn}(x)$ as
$$
\operatorname{sgn}(x)=
\begin{cases}
1, & \text{if } x>0,\\
-1, & \text{if } x<0,\\
0,  & \text{if } x=0.
\end{cases}
$$

Given an elliptic curve $E$ with traces of Frobenius $(a_2, a_3, \dots, a_{97})$, we define \[ \sgn(E)= (\sgn(a_2), \sgn(a_3), \dots, \sgn(a_{97})).\]

If we also specify a prime $p$, we define
 \[ \sgn(E, p)= (\sgn(a_2), \sgn(a_3), \dots, \widehat{\sgn(a_p)}, \dots, \sgn(a_{97})),\]
 where we exclude $\sgn(a_p)$ from the tuple $\sgn(E)$.

\subsection{Data} 
\label{subsec: data}

We use data derived from ECQ8 \cite{ECQ8}, which consists of isogeny classes of curves of conductor $N(E) < 10^8$. The experiments in Section \ref{sec: twists} use curves with conductor $N(E)  < 10^6$, and we refer to this dataset as ECQ6 \cite{ECQ6}. The specific data we use are the conductor $N(E)$, the first 25 traces of Frobenius $(a_2(E), a_3(E), \dots, a_{97}(E))$, and the equation of a curve in each isogeny class in order to compute the twist hash $h(E)$. In each experiment to predict $a_p(E)$,  we restrict ourselves to the set of curves of good reduction at $p$, i.e. the subset of curves whose conductor is not divisible by $p$.

The dataset generated in Section \ref{sec: beyond_twist} contains one curve per twist hash of elliptic curves with conductor up to $10^7$. This dataset includes the conductor $N(E)$, twist hash $h(E)$, and the list of the first 1229 traces of Frobenius $(a_p(E))_{p<10^4}$.

\section{Predicting $a_p(E)$}
\label{sec: twists}

In \cite{BCCHLLNP}, the authors trained  transformer-based models that could often correctly predict the magnitude of $a_p$ but struggled to distinguish the sign of $a_p$. Because quadratic twists share magnitudes $|a_p|$ at all but finitely many primes, this raises the question: could the high performance of these models be driven by the dataset allowing them to implicitly match curves to their quadratic twists, rather than learning deeper arithmetic properties?

To demonstrate that this dataset artifact alone is sufficient to generate successful predictions, we construct a baseline algorithm that predicts $a_p(E)$ and its sign by explicitly exploiting twist class redundancies. We evaluate this approach through two parallel experiments. In the first experiment, we use as input both the traces of Frobenius $a_p$ and the twist hash described in Section~\ref{subsec: ellcurves}. In the second experiment, we mimic the original experiments by only using the $a_p$ values to establish a proxy for the twist class. Specifically, the elliptic curves are instead grouped by the tuples \((|a_{p_1}|, \dots, |a_{p_k}|)\) consisting of the absolute values of traces at the last  $k$ primes. In both cases, we use curves from ECQ6.

\textit{Remark.} One could theoretically use $j$-invariants to partition the dataset into twist classes, since two quadratic twists $E$ and $E_\chi$ share the same $j$-invariant.  However, we would run into the issue that an isogeny class can, and often does, contain multiple isomorphism classes of curves over $\C$.

\subsection{Predicting $a_p(E)$ via Twist Hashes.}
\label{subsec: twists_exp1} In this experiment, we predict $a_p(E)$ by grouping the elliptic curves into twist classes using the twist hash $h(E)$. The steps are detailed in Algorithm \ref{alg:sign_matching_j}.

We use the ECQ6 dataset described in Section~\ref{subsec: data}. For each prime $p \in \{2, 3, 5, \dots, 97\}$, we restrict to curves of good reduction at $p$ (i.e., $p \nmid N(E)$), randomly shuffle the resulting set with a fixed seed, and reserve the last $10{,}000$ curves as the test set $T$ while the remaining curves form the training set $X$.

The input then consists of the training set $X$ and a test set $T$ of elliptic curves, all having good reduction at a prime $p < 100$. Each curve is represented by its twist hash $h(E)$ and its sequence of known traces of Frobenius $(a_q(E))_{q < 100,\, q \neq p}$. Using this data, the algorithm predicts the missing trace $a_p(E)$ for each test curve through the following steps.

\textbf{Step 0: Precomputation.} The algorithm first precomputes the probabilistic distribution $\mathcal{P}(p)$ of $a_p(E)$ on the training set, which serves as a fallback when no twist-hash match is available. Consequently, if a test curve's twist hash does not appear in the training set, $a_p(E)$ is randomly sampled from this distribution.

\begin{figure}[h]
\centering
\begin{tikzpicture}[
    every node/.style={font=\small},
    bubble/.style={circle, fill=red!80!black, draw=red!50!black,
                   inner sep=0pt, minimum size=8pt},
    lbl/.style={align=center, font=\small\bfseries},
]
    \draw[line width=0.6pt, draw=black!60, fill=black!8,
          rounded corners=4pt]
        (0,-0.12) rectangle (9,0.12);
    \node[bubble] (s0) at (2.25,0) {};
    \node[bubble] (s1) at (4.5,0)  {};
    \node[bubble] (s2) at (6.75,0) {};
    \node[lbl, above=10pt of s0] {Step 0:\\ Precomputation};
    \node[lbl, below=10pt of s1] {Step 1:\\ Sign-Pattern Indexing};
    \node[lbl, above=10pt of s2] {Step 2:\\ Matching Twist Hash\\ \& Majority Voting};
\end{tikzpicture}
\label{fig:pipeline}
\end{figure}

\textbf{Step 1: Sign-Pattern Indexing.} In this step we store the sign patterns between the curves in the training set that share same twist hash. For each pair of curves in $X$ sharing the same twist hash, it records the pointwise product of their sign vectors. The result is a lookup table $L$ that maps the sequence of relative signs at the known primes $q \neq p$ to the corresponding relative sign at the target prime $p$

\textit{Example.} To illustrate, consider two training curves $E_1, E_2 \in X$ belonging to the same twist class (i.e., $h(E_1) = h(E_2)$). For simplicity, suppose the algorithm is tasked to predict the third element in the sequence where $p=5$, i.e, $a_5$. Their traces of Frobenius differ only by sign. The value for $a_5$ is shown in bold:
\[
\left.\begin{aligned}
E_1 &= (-2,\; 0,\; \phantom{-}\mathbf{4},\; -3,\; 4,\, -5,\; \dots) \\
E_2 &= (\phantom{-}2,\; 1,\; \mathbf{-4},\; \phantom{-}3,\; 4,\, \phantom{-}5,\; \dots)
\end{aligned}\right\} \text{Twist curves}
\]

\textit{Step 1.1.} Extract the sign vectors:
\[
\begin{aligned}
\operatorname{sgn}(E_1) &= (-1,\; \phantom{-}0,\; \mathbf{+1},\; -1,\; +1,\, -1 \dots) \\
\operatorname{sgn}(E_2) &= (+1,\; +1,\; \mathbf{-1},\; +1,\; +1,\, +1 \dots)
\end{aligned}
\]

\textit{Step 1.2.} Drop the $a_5$ coordinate and take the pointwise product ($\ast$):
\[
\begin{array}{r@{\;}c@{\;}l}
\operatorname{sgn}(E_1, p)  &=& (-1,\; \phantom{-}0,\; -1,\; +1,\, -1 \dots) \\
\operatorname{sgn}(E_2, p)  &=& (+1,\; +1,\; +1,\; +1,\,+1\dots) \\
\hline
\operatorname{sgn}(E_1, E_2) &=& (-1,\; \phantom{-}0,\; -1,\; +1,\, -1\dots)
\end{array}
\]

\textit{Step 1.3.} Compute the target ($a_5$) sign relation:
\[
\operatorname{sgn}_5(E_1, E_2) \;=\; \operatorname{sgn}(a_5(E_1)) \cdot \operatorname{sgn}(a_5(E_2)) \;=\; (\mathbf{+1})(\mathbf{-1}) \;=\; -1.
\]

\textit{Step 1.4.} Append the feature-to-target pair to the lookup list:
\[
L \;\leftarrow\; \Big( \underbrace{(-1,\; 0,\; -1,\; +1,\ -1, \dots)}_{\operatorname{sgn}(E_1, E_2)},\;\; \underbrace{-1}_{\operatorname{sgn}_5(E_1, E_2)} \Big).
\]

\textbf{Step 2: Matching Twist Hash and Majority Voting.} For each test curve $E_T$, the algorithm retrieves all training curves $E_X$ that share its twist hash. For each matched $E_X$, the algorithm computes the relative sign pattern between $E_T$ and $E_X$ at the known primes. Querying the lookup table $L$ with this pattern yields the predicted relative sign at $p$, which in turn provides a candidate value for $a_p(E_T)$. The algorithm aggregates these candidates and selects the most frequent value as the final prediction for $a_p(E_T)$.

It is possible that a test curve $E_T$ matches a twist hash in the training set, but its specific relative sign pattern does not appear in $L$. In this case, the algorithm predicts $a_p(E_T)$ by sampling from the empirical distribution of $a_p(E_X)$ among those matched training curves. Finally, if multiple candidate values tie for the maximal frequency during the voting process, the algorithm breaks the tie by sampling from the tied candidates according to the initial probability distribution $\mathcal{P}(p)$.

\textit{Example.} To illustrate step 2, suppose a test curve $E_T \in T$ 
shares the same twist hash as $E_X \in X$, i.e., $h(E_T) = h(E_X)$:
\[
\begin{aligned}
E_T &= (\phantom{-}2,\; -1,\; \phantom{-}\mathbf{?},\; \phantom{-}5,\; 6,\; -3,\; \dots) \\
E_X &= (-2,\; \phantom{-}1,\; \mathbf{-4},\; -5,\; 0,\; \phantom{-}3,\; \dots)
\end{aligned}
\]

\textit{Step 2.1.}  Store $a_5(E_X)$ and compute the relative sign pattern $\operatorname{sgn}(E_T, E_X) = \operatorname{sgn}(E_T, 5) \ast \operatorname{sgn}(E_X, 5)$.
\[
\begin{array}{r@{\;}c@{\;}l}
\operatorname{sgn}(E_T, 5) &=& (+1,\; -1,\; +1,\; +1,\, -1,\, \dots) \\
\operatorname{sgn}(E_X, 5) &=& (-1,\; +1,\; -1,\; \phantom{-}0,\, +1,\, \dots) \\
\hline
\operatorname{sgn}(E_T, E_X) &=& (-1,\; -1,\; -1,\; \phantom{-}0,\, -1,\, \dots)
\end{array}
\]

\textit{Step 2.2.} Query this relative sign pattern $\operatorname{sgn}(E_T, E_X)$ in the lookup table $L$. This retrieves all pairs $(E_i, E_j)$ such that
$\operatorname{sgn}(E_T, E_X) = \operatorname{sgn}(E_i, E_j)$ and stores the occurrences of $\operatorname{sgn}_5(E_i, E_j)$ into a multi-set $M$. Suppose there are $100$ such matches, distributed as follows:
$$
M = \{\underbrace{-1, \dots, -1}_{80 \text{ times}}, \underbrace{0, \dots, 0}_{20 \text{ times}}\}
$$

\textit{Step 2.3.} Count the occurrences in the multi-set $M$ to form a frequency dictionary $freq$. 
$$
freq \;\leftarrow\; \{-1: 80,\; 0: 20\}
$$

\textit{Step 2.4.} Compute the candidate values for the missing trace $a_5(E_T)$. Each key $s \in freq$ encodes a relative sign, not an absolute trace. To recover candidate values for $a_5(E_T)$, we multiply each key in our frequency dictionary by the known value $a_5(E_X)=-4$.

\[
\begin{aligned}
freq &\;\leftarrow\; \{-1 * a_5(E_X): 80,\; 0 * a_5(E_X): 20\} \\
\text{candidates} &= \{4: 80,\; 0: 20\} \\
\end{aligned}
\]

\textit{Step 2.5.} Repeat \textit{Steps 2.1–2.5} for all training curves $E_X \in X$ sharing the twist hash $h(E_T)$, and accumulate the resulting candidate frequencies into a global map. The final predicted trace, $\widehat{a_5}(E_T)$ is the candidate with the maximal aggregate vote count:
\[
a_5(E_T) \;=\; \arg\max_{n \in \mathbb{Z}} freq[n] \;=\; 4.
\]

\begin{figure}[htbp!]
    \centering
    \includegraphics[width=\textwidth]{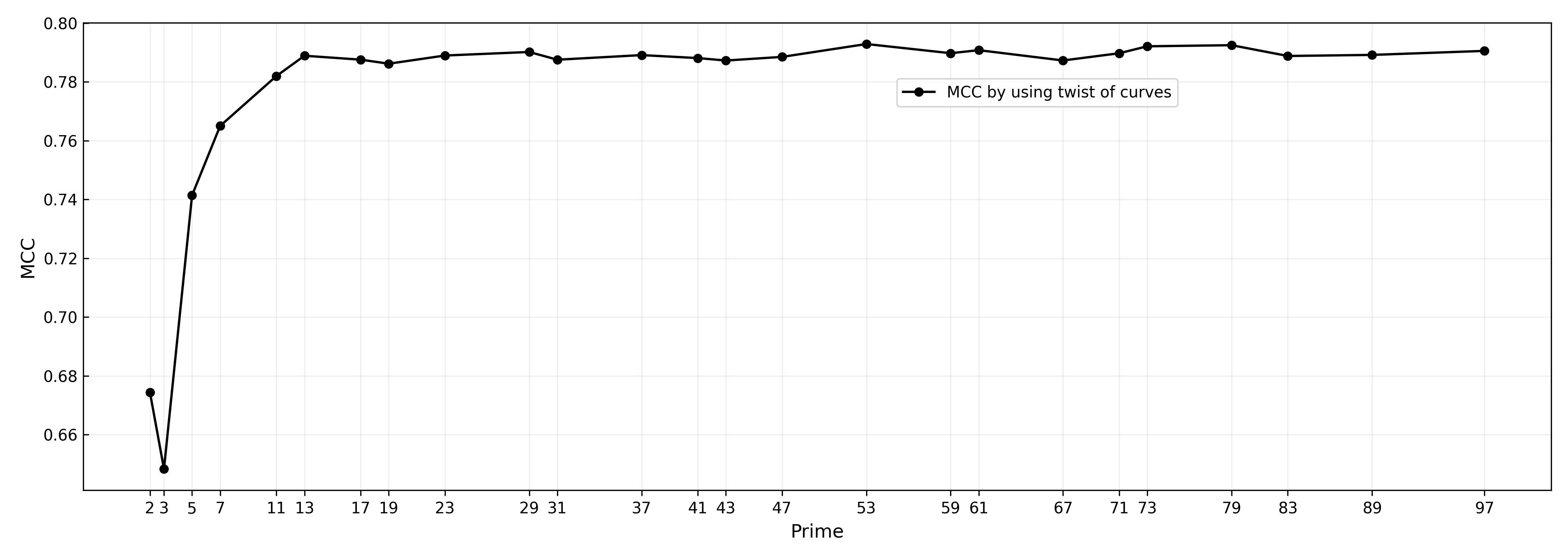}
    \caption{MCC of the sign-based twist class matching algorithm (Algorithm~\ref{alg:sign_matching_j}) 
for predicting $a_p(E)$ across primes $p \in \{2,3, \dots,  97\}$, using the exact twist hash to group 
curves.}
    \label{fig:mcc_twist}
\end{figure}

Figure~\ref{fig:mcc_twist} shows the results of Experiment 1, where  Algorithm~\ref{alg:sign_matching_j} starts with an MCC of $0.6744$ when predicting $a_2(E)$ and $0.6483$ when predicting $a_3(E)$, and plateaus at $\approx 0.79$ for larger primes starting at $p=11$. We specifically highlight the performance at $a_{97}(E)$ (MCC = $0.7905$) to provide a direct comparison with the transformer models of \cite{BCCHLLNP}, which were evaluated on $a_2(E)$, $a_3(E)$, and $a_{97}(E)$. Those models achieved MCCs of $0.5822$, $0.5205$, and $0.5137$, respectively. Thus, our explicit twist-matching approach significantly outperforms the previously trained transformer models across the board.

It is possible that the transformer models primarily derived their predictive capabilities by implicitly recognizing these twist classes. However, because the explicit twist hash was not provided as an input feature in \cite{BCCHLLNP}, we cannot observe this mechanism directly. This motivates Experiment 2, where instead of using the twist hash directly,  we group curves by tuples of absolute values $(|a_{p_1}(E)|, \dots, |a_{p_k}(E)| )$ as a proxy 
for the twist class. We then examine how closely this proxy recovers the performance  seen in Experiment 1.

\subsection{Approximating the twist hash partition from $(|a_{p_1}|, \dots, |a_{p_k}|)$}
We seek to closely approximate the exact twist class grouping. To achieve this, we use the absolute values of the traces of Frobenius $(|a_{p_1}(E)|, \dots, |a_{p_k}(E)|)$ at the largest $k$ primes as a proxy. We rely on larger primes because they are less likely to divide the conductor of the curve, meaning the curves more consistently exhibit good reduction at these primes. This makes the magnitudes of traces of Frobenius at larger primes a robust and reliable indicator of the twist class. If this heuristic is correct, a partition built from the largest primes should much more closely mirror the exact twist-hash partition than one built from the smallest primes.

\begin{figure}[htbp!]
    \centering
    \includegraphics[width=\textwidth]{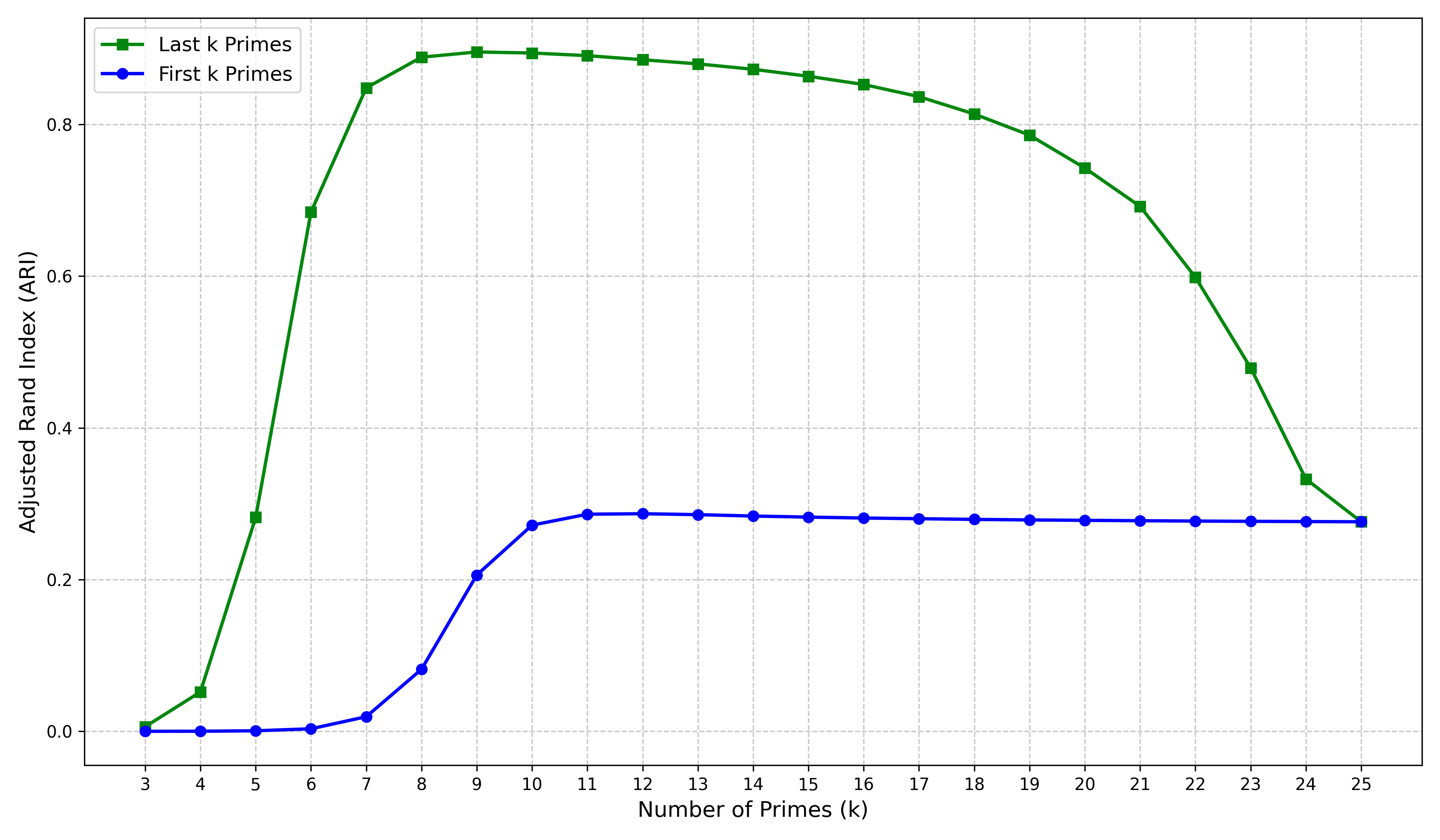}
    \caption{Adjusted Rand Index (ARI) between the twist-hash partition and the partition of $\mathrm{ECQ6}$ induced by tuples of absolute values $(|a_{p_1}(E)|, \dots, |a_{p_k}(E)|)$ as a function of $k$. The green curves uses the $k$ largest primes below $100$ while the blue uses the $k$ smallest. The last $k$ partition recovers the twist-hash classes substantially better, with ARI peaking near $0.85$ on $k \in [7, 16]$ while the first-$k$ partition saturates near $0.28$}
    \label{fig:last_first_k_primes}
\end{figure}

We therefore consider two  choices: the $k$ smallest and the $k$ largest primes below 100. Figure~\ref{fig:last_first_k_primes} compares the resulting partitions against the twist-hash partition. The $k$ largest primes recover the twist-hash partition substantially better, with ARI peaking near 0.85. Conversely, using the first $k$ primes yields rapidly diminishing results. This stark contrast validates our hypothesis. 

In our setting, ARI is the most informative overall score, since it continues to distinguish between good class recovery and later overfragmentation as \(k\) grows. By contrast, homogeneity, completeness, and V-measure are useful supporting summaries, but they can saturate rapidly in the presence of many small or singleton twist-hash classes.  Furthermore, special care is needed because the true twist-hash partition contains many singleton classes. On this fully atomized singleton subset (where each curve forms its own class), the ARI behaves almost like a collision detector for the trace representation. Therefore, the genuine shared-class recovery problem is most accurately reflected by the ARI evaluated strictly on the multi-instance subset (i.e., the subset of curves whose true twist-hash class has a size of at least two).

\begin{figure}[htbp!]
    \centering
    \includegraphics[width=\textwidth]{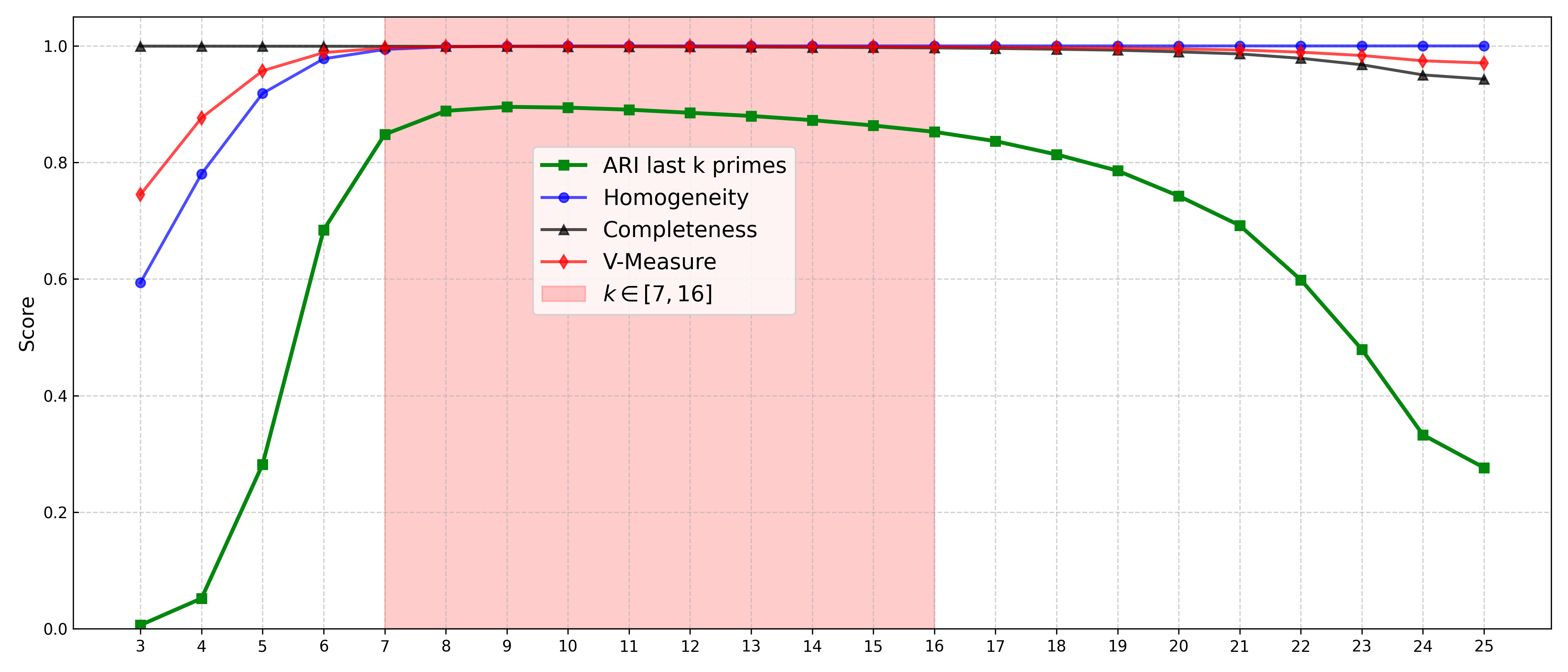}
    \caption{ARI, Homegeneity, Completeness and V-measure of the partitioning using the tuples of absolute values $(|a_{p_1}(E)|, \dots, |a_{p_k}(E)|)$. The range of $k \in [7, 16]$ yields the maximum score for all these metrics.}
    \label{fig:all_stats_clustring}
\end{figure}

Figure~\ref{fig:all_stats_clustring} reports the homogeneity, completeness, and V-measure for the largest $k$ primes, alongside the ARI. The three entropy-based scores rise quickly with $k$ and are essentially saturated on $k \in [7, 16]$, while the ARI is also near its maximum on the same range. We therefore restrict attention to $k \in [7, 16]$ for our proxy models in the remainder of the section.

\subsection{A proxy algorithm for predicting $a_p(E)$} 
To construct a proxy algorithm for the twist-hash model of Section~\ref{subsec: twists_exp1}, we approximate the twist-hash partition using the tuples of absolute values $(|a_{p_1}(E)|, \dots, |a_{p_k}(E)|)$, and use these tuples as proxy keys in place of the twist hash in Algorithm~\ref{alg:sign_matching_j}. The resulting procedure is recorded as Algorithm~\ref{alg:sign_matching_proxy} in the Appendix.

\begin{figure}[htbp!]
    \centering
    \includegraphics[width=\textwidth]{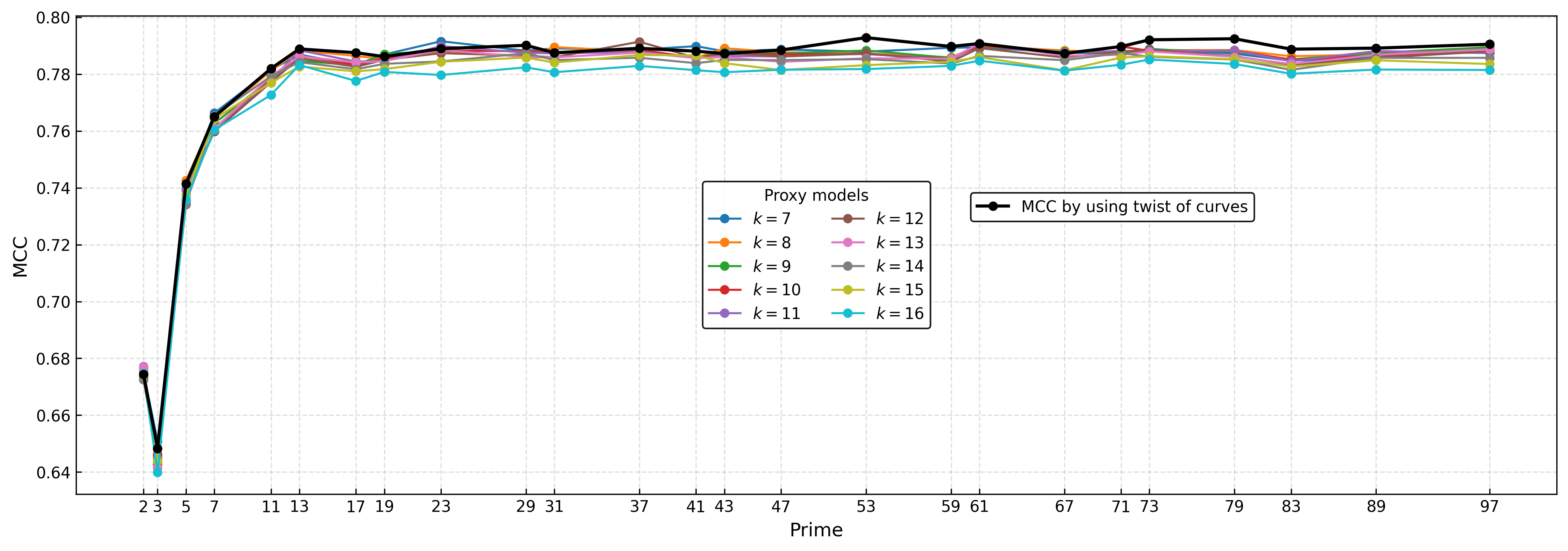}
    \caption{Comparison of proxy models to the twist-hash model.}
    \label{fig:proxymodels}
\end{figure}

Algorithm~\ref{alg:sign_matching_proxy} differs from Algorithm~\ref{alg:sign_matching_j} only in the key used to group curves. In place of the twist hash $h(E)$, it uses the proxy key $k(E) = \big(|a_{p_1}(E)|, \dots, |a_{p_k}(E)|\big)$ formed from the absolute values of the traces at the $k$ largest primes below 100. Steps 0, 1 and 2, as listed in Section~\ref{subsec: twists_exp1}, are otherwise unchanged, with proxy keys taking the role of twist hashes throughout.

We run Algorithm~\ref{alg:sign_matching_proxy} for each $k\in[7,16]$ on ECQ6, using the same train/test split as in Experiment~1. Figure~\ref{fig:proxymodels} shows that the per-prime MCC of every proxy model closely tracks that of the twist-hash model across all primes $p\in[2,97]$. To identify the proxy closest to the twist-hash model, we compute the MAE and RMSE between their per-prime MCC values. Table~\ref{tab:proxy_proximity} reports these statistics, and $k=8$ attains the smallest MAE.

\begin{table}[htbp!]
\centering
\begin{tabular}{ccccc}
\hline
\textbf{Proxy (k)} & \textbf{MAE} & \textbf{RMSE} & \textbf{Max Difference} & \textbf{Mean Difference} \\
\hline
7  & 0.0020 & 0.0025 & 0.0053 & -0.0010 \\
8  & 0.0019 & 0.0024 & 0.0056 & -0.0013 \\
9  & 0.0022 & 0.0026 & 0.0046 & -0.0020 \\
10 & 0.0023 & 0.0029 & 0.0058 & -0.0021 \\
11 & 0.0025 & 0.0030 & 0.0057 & -0.0022 \\
12 & 0.0032 & 0.0037 & 0.0062 & -0.0030 \\
13 & 0.0032 & 0.0038 & 0.0073 & -0.0030 \\
14 & 0.0041 & 0.0045 & 0.0076 & -0.0041 \\
15 & 0.0047 & 0.0052 & 0.0097 & -0.0047 \\
16 & 0.0071 & 0.0074 & 0.0111 & -0.0070 \\
\hline
\end{tabular}
\caption{Proximity of proxy models to the twist-hash model with MAE (Mean Absolute Error) and RMSE (Root Mean Squared Error). The closest proxy model by MAE is $k=8$.}
\label{tab:proxy_proximity}
\end{table}

Finally, we evaluate how frequently the $k=8$ proxy model and the twist-hash model yield identical predictions for individual test curves. Measuring this disagreement is crucial, because even if two models achieve similar aggregate performance metrics, a high disagreement rate would indicate that their predictive power stems from different sources, causing them to predict entirely different subsets of curves correctly. As summarized in Figure~\ref{fig:proxy_disagreement} and Table~\ref{tab:disagreement_k8} in the Appendix, the disagreement rate between the two models exceeds $10\%$ at $p=2$, but drops below $3\%$ for $p \ge 29$ and plateaus at $\approx2\%$ for higher primes. This strong, instance-level alignment confirms that our preferred choice of $k=8$ is a highly faithful proxy for the explicit twist class.

\begin{figure}[htbp!]
    \centering
    \includegraphics[width=\textwidth]{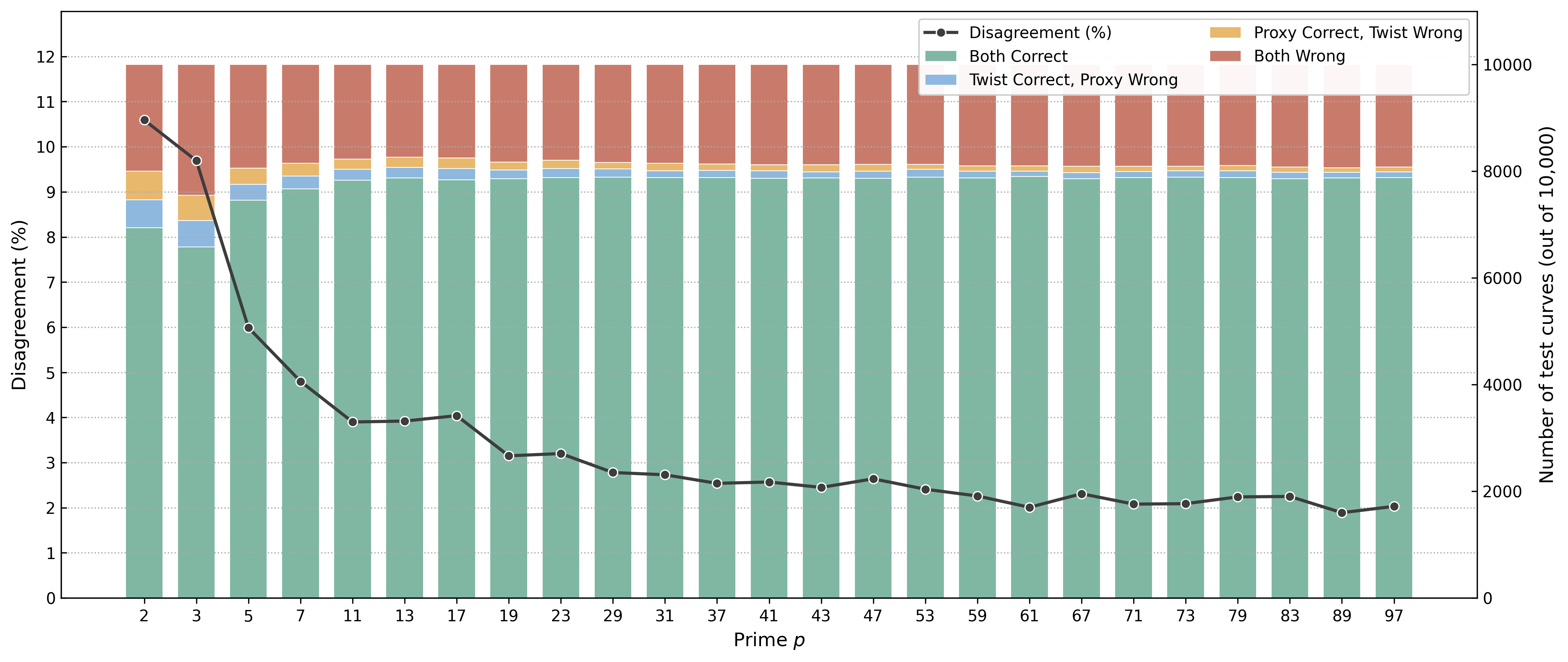}
    \caption{Agreement and disagreement between the twist model and proxy model $(k=8)$}
    \label{fig:proxy_disagreement}
\end{figure}

\section{Towards Generalization Beyond Twist Classes}
\label{sec: beyond_twist}

\subsection{Evaluation of trained transformer models on unseen twist classes}
\label{subsec: unseen}
The previous sections demonstrate that twist classes provide a massive predictive advantage. However, this leaves open the question of whether the transformer models' predictions were based on the presence of multiple curves from the same twist class in both training and test set. In this subsection we investigate the models to predict $a_p(E)$ for $p=2,3$ and $97$ trained in \cite{BCCHLLNP} specifically on unseen twist classes.

For each prime, we compared the observed exact-value MCC against a training-marginal null baseline. Under this null hypothesis, the baseline guesses are generated by sampling from the probability distribution $\mathcal{P}(p)$ of the training data. Evaluating at a significance level of $\alpha=0.01$, the model demonstrated statistically significant exact-value predictive power across all three primes ($p < 0.001$). This proves that the transformer models capture structural signal beyond random prediction from the empirical distribution.

In \cite{BCCHLLNP}, the authors noticed that reducing the prediction modulo $\ell$ for $\ell=2, 4$ resulted in a high performance of such models. We reproduce the experiment for the raw value, now reducing modulo $2$, $3$ and $4$.  For all three primes, the model exhibited its strongest predictive performance on 2-adic features, with parity ($a_p \bmod 2$) consistently yielding the highest performance (MCCs between $0.293$ and $0.381$), followed closely by the modulo 4 task, all being statistically significant ($p<0.001$). While predictions reduced modulo 3 remained statistically significant across the board ($p<0.001$), the model's accuracy on 3-adic features was persistently weaker than its performance modulo $2$ or modulo $4$. The experiment details can be found in Table \ref{tab:multi_prime_results}.

\begin{table}[hbt!]
\centering
\begin{tabular}{llrrl}
\toprule
\textbf{Prime} & \textbf{Prediction Task} & \textbf{Observed MCC} & \textbf{Null Mean} & \textbf{$p$-value} \\
\midrule
$p=2$ ($n=2904$)
& Exact Value & $0.102$ & $0.0000$ & $p < 0.001^\ast$ \\
& Parity ($a_p \bmod 2$) & $0.303$ & $-0.0005$ & $p < 0.001^\ast$ \\
& Modulo 3 & $0.065$ & $-0.0001$ & $p < 0.001^\ast$ \\
& Modulo 4 & $0.108$ & $0.0009$ & $p < 0.001^\ast$ \\
\midrule
$p=3$ ($n=3293$)
& Exact Value & $0.117$ & $-0.0005$ & $p < 0.001^\ast$ \\
& Parity ($a_p \bmod 2$) & $0.381$ & $0.0000$ & $p < 0.001^\ast$ \\
& Modulo 3 & $0.059$ & $-0.0001$ & $p < 0.001^\ast$ \\
& Modulo 4 & $0.183$ & $0.0001$ & $p < 0.001^\ast$ \\
\midrule
$p=97$ ($n=2031$)
& Exact Value & $0.023$ & $0.0001$ & $p < 0.001^\ast$ \\
& Parity ($a_p \bmod 2$) & $0.294$ & $0.0001$ & $p < 0.001^\ast$ \\
& Modulo 3 & $0.047$ & $-0.0005$ & $p \approx 0.004^\ast$ \\
& Modulo 4 & $0.138$ & $-0.0003$ & $p < 0.001^\ast$ \\
\bottomrule
\multicolumn{5}{l}{\footnotesize $^\ast$Statistically significant at the $\alpha = 0.01$ level.}
\end{tabular}
\caption{Statistical significance of transformer predictions across multiple primes ($p=2, 3, 97$) on unseen twist classes, evaluated against a training-marginal null baseline ($1{,}000$ draws weighted by training set distributions).}
\label{tab:multi_prime_results}
\end{table}
\subsection{A dataset with a unique representative per twist class}

As seen in Section \ref{sec: twists}, models that utilize the absolute values of the traces of Frobenius of the largest $k$ primes as a proxy for the quadratic twist class of an elliptic curve, and use this feature as a look-up key, significantly outperform any trained models from \cite{BCCHLLNP}. Moreover, as described in subsection \ref{subsec: unseen}, the ML models performed substantially worse on the slice of the test set of curves with no quadratic twists in the training set. Despite this, the transformer models seem to learn additional properties of the sequences of traces of Frobenius. To account for the information contained in the presence of multiple twists, in this section we introduce a dataset, derived from ECQ8, consisting of exactly one representative per quadratic twist class of curves with conductor $N(E)<10^7$.

In particular, each curve in ECQ8 is represented by a choice of coefficients $a_1, a_2, a_3, a_4$ and $a_6$ describing the Weierstrass equation $y^2+a_1xy+a_3y=x^3+a_2x^2+a_4x+a_6$ of a curve in this isogeny class. Note that any two curves in the same isogeny class have the same twist hash. Once we compute the twist hash $h(E)$, we only keep the isogeny class showing up first, which are ordered in $\ECQ8$ in ascending order by conductor. Finally, we record the following data of isogeny invariants, consisting of the conductor $N(E)$, the twist hash $h(E)$, and the first 1229 traces of Frobenius $(a_p(E))_{p<10^4}$. The dataset \cite{ECQ7} is publicly available.

\bibliographystyle{plain}
\bibliography{traces}

@inproceedings{RosenbergHirschberg2007,
  author    = {Rosenberg, Andrew and Hirschberg, Julia},
  title     = {V-Measure: A Conditional Entropy-Based External Cluster Evaluation Measure},
  booktitle = {Proceedings of the 2007 Joint Conference on Empirical Methods in Natural Language Processing and Computational Natural Language Learning},
  pages     = {410--420},
  year      = {2007}
}

@article{gorodkin2004comparing,
  title={Comparing two K-category assignments by a K-category correlation coefficient},
  author={Gorodkin, Jan},
  journal={Computational biology and chemistry},
  volume={28},
  number={5-6},
  pages={367--374},
  year={2004},
  publisher={Elsevier}
}

@article{matthews1975comparison,
  title={Comparison of the predicted and observed secondary structure of T4 phage lysozyme},
  author={Matthews, Brian W},
  journal={Biochimica et Biophysica Acta (BBA)-Protein Structure},
  volume={405},
  number={2},
  pages={442--451},
  year={1975},
  publisher={Elsevier}
}

@book {silvermanArithmetic,
    AUTHOR = {Silverman, Joseph H.},
     TITLE = {The arithmetic of elliptic curves},
    SERIES = {Graduate Texts in Mathematics},
    VOLUME = {106},
   EDITION = {Second},
 PUBLISHER = {Springer, Dordrecht},
      YEAR = {2009},
     PAGES = {xx+513},
      ISBN = {978-0-387-09493-9},
   MRCLASS = {11-02 (11G05 11G20 14H52 14K15)},
  MRNUMBER = {2514094},
MRREVIEWER = {Vasil\cprime\ \=I.\ Andr\=\i\u ichuk},
       DOI = {10.1007/978-0-387-09494-6},
       URL = {https://doi.org/10.1007/978-0-387-09494-6},
}

@dataset{ECQ6,
  SHORTHAND = {ECQ6},
  author    = {Costa, Edgar},
  title     = {{Frobenious traces for a set of isogeny classes of elliptic curves of conductor up to $10^6$}},
  publisher = {Zenodo},
  year      = {2025},
  month     = jun,
  version   = {v1},
  doi       = {10.5281/zenodo.15777474},
  url       = {https://zenodo.org/records/15777475}
}

@dataset{ECQ8,
  SHORTHAND = {ECQ8},
  author    = {Sutherland, Andrew Victor},
  title     = {{A set of isogeny classes of elliptic curves of conductor up to $10^8$}},
  publisher = {Zenodo},
  year      = {2024},
  month     = sep,
  version   = {v1},
  doi       = {10.5281/zenodo.14847809},
  url       = {https://doi.org/10.5281/zenodo.14847809}
}

@article{BCCHLLNP,
     AUTHOR = {Babei, Angelica and Charton, Fran{\c{c}}ois and Costa, Edgar and Huang, Xiaoyu and Lee, Kyu-Hwan and Lowry-Duda, David and Narayanan, Ashvni and Pozdnyakov, Alexey},
      TITLE = {Learning Euler factors of elliptic curves},
    JOURNAL = {Advances in Theoretical and Mathematical Physics},
      PAGES = {2327--2351},
       YEAR = {2025},
  VOLUME = {29},
  NUMBER={8},
        DOI = {10.4310/ATMP.260412180132},
}

@article{BBSHS24,
     AUTHOR = {Babei, Angelica and Banwait, Barinder S. and Fing, AJ and Huang, Xiaoyu and Singh, Deependra},
      TITLE = {Machine learning approaches to the Shafarevich-Tate group of elliptic curves},
    JOURNAL = {Advances in Theoretical and Mathematical Physics},
      PAGES = {pp. 2353-2379},
       YEAR = {2025},
  VOLUME = {29},
  NUMBER={8},
        DOI = {10.4310/ATMP.260412210820},
}

@article {Magma,
    AUTHOR = {Bosma, Wieb and Cannon, John and Playoust, Catherine},
     TITLE = {The {M}agma algebra system. {I}. {T}he user language},
      NOTE = {Computational algebra and number theory (London, 1993)},
   JOURNAL = {J. Symbolic Comput.},
  FJOURNAL = {Journal of Symbolic Computation},
    VOLUME = {24},
      YEAR = {1997},
    NUMBER = {3-4},
     PAGES = {235--265},
      ISSN = {0747-7171},
   MRCLASS = {68Q40},
  MRNUMBER = {MR1484478},
       DOI = {10.1006/jsco.1996.0125},
       URL = {http://dx.doi.org/10.1006/jsco.1996.0125},
}

@misc{SutherlandMagma,
  author = {Sutherland, Andrew V.},
  title = {Magma repository},
  year = {2023},
  publisher = {GitHub},
  journal = {GitHub repository},
  howpublished = {\url{https://github.com/AndrewVSutherland/Magma}},
}

@article{BSSVY, 
title={A database of genus-2 curves over the rational numbers}, 
volume={19}, 
DOI={10.1112/S146115701600019X}, 
number={A}, 
journal={LMS Journal of Computation and Mathematics}, 
author={Booker, Andrew R. and Sijsling, Jeroen and Sutherland, Andrew V. and Voight, John and Yasaki, Dan}, 
year={2016}, 
pages={235–254}
}

@dataset{ECQ7,
  author       = {Babei, Angelica},
  title        = {Frobenius traces for a set of (quadratic) twist
                   classes of elliptic curves of conductor up to $10^7$
                  },
  month        = may,
  year         = 2026,
  publisher    = {Zenodo},
  version      = {v0},
  doi          = {10.5281/zenodo.20129317},
  url          = {https://doi.org/10.5281/zenodo.20129317},
}

@article{HubertArabie1985,
  author  = {Hubert, Lawrence and Arabie, Phipps},
  title   = {Comparing Partitions},
  journal = {Journal of Classification},
  volume  = {2},
  number  = {1},
  pages   = {193--218},
  year    = {1985}
}

@misc{lmfdb,
  shorthand    = {LMFDB},
  author       = {The {LMFDB Collaboration}},
  title        = {The {L}-functions and modular forms database},
  howpublished = {\url{https://www.lmfdb.org}},
  year         = {2024},
  note         = {[Online; accessed 29 December 2024]},
}

@misc{BBLLD,
  title={Murmurations of modular forms in the weight aspect},
  author={Jonathan Bober and Andrew R. Booker and Min Lee and David Lowry-Duda},
  year={2023},
  howpublished="\url{http://arxiv.org/abs/2310.07746v1}",
  note={{arXiv:math.NT:2310.07746v1}},
}

@misc{BLLD+,
  title={Murmurations of {M}aass forms},
  author={Andrew R. Booker and Min Lee and David Lowry-Duda and Andrei Seymour-Howell and Nina Zubrilina},
  year={2024},
  howpublished="\url{http://arxiv.org/abs/2409.00765v1}",
  note={{arXiv:math.NT:2409.00765v1}},
}

@misc{BKM,
  title={Murmurations of {M}estre-{N}agao sums},
  author={Zvonimir Bujanović and Matija Kazalicki and Lukas Novak},
  year={2024},
  howpublished="\url{http://arxiv.org/abs/2403.17626v1}",
  note={{arXiv:math.NT:2403.17626v1}},
}

@misc{Cowan,
  title={Murmurations and explicit formulas},
  author={Alex Cowan},
  year={2023},
  howpublished="\url{http://arxiv.org/abs/2306.10425v2}",
  note={{arXiv:math.NT:2306.10425v2}},
}

@misc{C1,
  title={Murmurations and ratios conjectures},
  author={Alex Cowan},
  year={2024},
  howpublished="\url{http://arxiv.org/abs/2408.12723v1}",
  note={{arXiv:math.NT:2408.12723v1}},
}

@article{LOP25,
  title={Murmurations of {D}irichlet characters},
  author={Lee, Kyu-Hwan and Oliver, Thomas and Pozdnyakov, Alexey},
  journal={International Mathematics Research Notices},
  volume={2025},
  number={1},
  year={2025},
  publisher={Oxford University Press}
}

@article{SawinSutherland2025,
  author  = {Sawin, Will and Sutherland, Andrew V.},
  title   = {Murmurations for Elliptic Curves Ordered by Height},
  journal = {arXiv preprint},
  eprint  = {2504.12295},
  archivePrefix = {arXiv},
  primaryClass  = {math.NT},
  year    = {2025},
}

@misc{Martin,
  title={Distribution of local signs of modular forms and murmurations of {F}ourier
  coefficients},
  author={Kimball Martin},
  year={2024},
  howpublished="\url{http://arxiv.org/abs/2409.02338v1}",
  note={{arXiv:math.NT:2409.02338v1}},
}

@article{HLOP24,
     AUTHOR = {Yang-Hui He and Kyu-Hwan Lee and Thomas Oliver and Alexey Pozdnyakov},
      TITLE = {Murmurations of Elliptic Curves},
    JOURNAL = {Experimental Mathematics},
      PAGES = {1--13},
       YEAR = {2024},
  PUBLISHER = {Taylor \& Francis},
        DOI = {10.1080/10586458.2024.2382361},
}

@misc{Zubrilina,
  title={Murmurations},
  author={Nina Zubrilina},
  year={2023},
  howpublished="\url{http://arxiv.org/abs/2310.07681v1}",
  note={{arXiv:math.NT:2310.07681v1}},
}

@article {KV23,
    AUTHOR = {Kazalicki, Matija and Vlah, Domagoj},
     TITLE = {Ranks of elliptic curves and deep neural networks},
   JOURNAL = {Res. Number Theory},
  FJOURNAL = {Research in Number Theory},
    VOLUME = {9},
      YEAR = {2023},
    NUMBER = {3},
     PAGES = {Paper No. 53, 21},
      ISSN = {2522-0160,2363-9555},
   MRCLASS = {11G05 (68T07)},
  MRNUMBER = {4611298},
MRREVIEWER = {John\ M.\ Voight},
       DOI = {10.1007/s40993-023-00462-w},
       URL = {https://doi.org/10.1007/s40993-023-00462-w},
}

@misc{AHLOS,
  AUTHOR = {M. Amir and Y.-H.~He and K.-H.~Lee and T.~Oliver and E. Sultanow},
  TITLE = {Machine Learning Class Numbers of Real Quadratic Fields},
  year={2022},
  howpublished="\url{https://arxiv.org/pdf/2209.09283}",
  note={{arXiv:math.NT:2209.09283}},
}

@misc{KS25,
  AUTHOR = {K.-H.~Lee and S. Lee},
  TITLE = {Machines Learn Number Fields, But How? The Case of Galois Groups},
  year={2025},
  howpublished="\url{https://arxiv.org/pdf/2508.06670}",
  note={{arXiv:math.NT:2508.06670}},
}

\newpage

\section*{Appendix}

\begin{table}[htbp!]
\centering
\caption{Comparison of Twist Hash and Proxy ($k=8$) Models}
\label{tab:merged_metrics}
\resizebox{\textwidth}{!}{%
\begin{tabular}{cc|cccc|cccc}
\hline
\multirow{2}{*}{\textbf{Prime}} & 
\multirow{2}{*}{\textbf{\makecell{\# of Good \\ Conductors}}} & 
\multicolumn{4}{c|}{\textbf{Twist Hash Model}} & 
\multicolumn{4}{c}{\textbf{Proxy Model ($k=8$)}} \\
\cline{3-10}
& & 
\textbf{\makecell{Probabilistically \\ Predicted$^*$}} & 
\textbf{\makecell{Overall \\ Correct$^*$}} & 
\textbf{\makecell{Deterministic \\ MCC}} & 
\textbf{\makecell{Overall \\ MCC}} & 
\textbf{\makecell{Probabilistically \\ Predicted$^*$}} & 
\textbf{\makecell{Overall \\ Correct$^*$}} & 
\textbf{\makecell{Deterministic \\ MCC}} & 
\textbf{\makecell{Overall \\ MCC}} \\
\hline
2  & 722205  & 3177 & 7466 & 0.9834 & 0.6744 & 3198 & 7481 & 0.9830 & 0.6762 \\
3  & 1210417 & 3489 & 7079 & 0.9839 & 0.6483 & 3513 & 7058 & 0.9847 & 0.6458 \\
5  & 1996257 & 2561 & 7754 & 0.9909 & 0.7414 & 2579 & 7765 & 0.9910 & 0.7426 \\
7  & 2466121 & 2328 & 7911 & 0.9907 & 0.7650 & 2344 & 7909 & 0.9909 & 0.7646 \\
11 & 2932416 & 2188 & 8035 & 0.9935 & 0.7819 & 2202 & 8027 & 0.9931 & 0.7811 \\
13 & 3057496 & 2165 & 8073 & 0.9979 & 0.7889 & 2176 & 8069 & 0.9979 & 0.7884 \\
17 & 3207083 & 2171 & 8050 & 0.9967 & 0.7875 & 2185 & 8040 & 0.9965 & 0.7864 \\
19 & 3259866 & 2140 & 8021 & 0.9952 & 0.7862 & 2157 & 8014 & 0.9953 & 0.7853 \\
23 & 3330210 & 2124 & 8049 & 0.9949 & 0.7889 & 2139 & 8043 & 0.9948 & 0.7882 \\
29 & 3396401 & 2133 & 8044 & 0.9970 & 0.7901 & 2145 & 8008 & 0.9969 & 0.7862 \\
31 & 3412665 & 2128 & 8007 & 0.9950 & 0.7875 & 2141 & 8026 & 0.9957 & 0.7895 \\
37 & 3447726 & 2144 & 8017 & 0.9978 & 0.7891 & 2157 & 8007 & 0.9982 & 0.7880 \\
41 & 3465643 & 2142 & 8009 & 0.9984 & 0.7881 & 2156 & 7986 & 0.9980 & 0.7856 \\
43 & 3472612 & 2144 & 7992 & 0.9976 & 0.7872 & 2160 & 8009 & 0.9981 & 0.7891 \\
47 & 3485303 & 2138 & 8003 & 0.9968 & 0.7885 & 2152 & 7995 & 0.9966 & 0.7877 \\
53 & 3499565 & 2145 & 8038 & 0.9995 & 0.7928 & 2154 & 7985 & 0.9988 & 0.7872 \\
59 & 3511808 & 2132 & 8006 & 0.9977 & 0.7897 & 2146 & 7970 & 0.9976 & 0.7859 \\
61 & 3515014 & 2138 & 8006 & 0.9984 & 0.7908 & 2155 & 7999 & 0.9987 & 0.7900 \\
67 & 3524379 & 2160 & 7974 & 0.9988 & 0.7872 & 2173 & 7985 & 0.9992 & 0.7883 \\
71 & 3529542 & 2119 & 7999 & 0.9963 & 0.7897 & 2136 & 7979 & 0.9957 & 0.7876 \\
73 & 3531417 & 2133 & 8011 & 0.9991 & 0.7921 & 2149 & 7976 & 0.9992 & 0.7884 \\
79 & 3536694 & 2139 & 8013 & 0.9976 & 0.7924 & 2155 & 7975 & 0.9976 & 0.7885 \\
83 & 3540455 & 2135 & 7983 & 0.9976 & 0.7888 & 2151 & 7960 & 0.9969 & 0.7863 \\
89 & 3544529 & 2143 & 7983 & 0.9992 & 0.7891 & 2163 & 7962 & 0.9992 & 0.7869 \\
97 & 3550728 & 2141 & 7987 & 0.9984 & 0.7905 & 2159 & 7978 & 0.9991 & 0.7896 \\
\hline
\end{tabular}%
} 
\vspace{1.5ex} 
\raggedright
\textit{\footnotesize $^*$Denotes the number of correct predictions out of 10,000 total test samples.}
\end{table}

\begin{table}[htbp!]
\centering
\footnotesize 
\setlength{\tabcolsep}{4pt}
\caption{Agreement and disagreement between the twist model and proxy model (k=8).}
\label{tab:disagreement_k8}
\begin{tabular}{cccccc}
\hline
\textbf{Prime} & \textbf{\makecell{Both \\ Correct}} & \textbf{\makecell{Twist Correct \\ Proxy Wrong}} & \textbf{\makecell{Proxy Correct \\ Twist Wrong}} & \textbf{\makecell{Both \\ Wrong}} & \textbf{\makecell{Disagreement \\ (\%)}} \\
\hline
2  & 6944 & 522 & 537 & 1997 & 10.59 \\
3  & 6584 & 495 & 474 & 2447 &  9.69 \\
5  & 7460 & 294 & 305 & 1941 &  5.99 \\
7  & 7670 & 241 & 239 & 1850 &  4.80 \\
11 & 7836 & 199 & 191 & 1774 &  3.90 \\
13 & 7875 & 198 & 194 & 1733 &  3.92 \\
17 & 7843 & 207 & 197 & 1753 &  4.04 \\
19 & 7860 & 161 & 154 & 1825 &  3.15 \\
23 & 7886 & 163 & 157 & 1794 &  3.20 \\
29 & 7887 & 157 & 121 & 1835 &  2.78 \\
31 & 7880 & 127 & 146 & 1847 &  2.73 \\
37 & 7885 & 132 & 122 & 1861 &  2.54 \\
41 & 7869 & 140 & 117 & 1874 &  2.57 \\
43 & 7878 & 114 & 131 & 1877 &  2.45 \\
47 & 7867 & 136 & 128 & 1869 &  2.64 \\
53 & 7891 & 147 &  94 & 1868 &  2.41 \\
59 & 7875 & 131 &  95 & 1899 &  2.26 \\
61 & 7902 & 104 &  97 & 1897 &  2.01 \\
67 & 7864 & 110 & 121 & 1905 &  2.31 \\
71 & 7885 & 114 &  94 & 1907 &  2.08 \\
73 & 7889 & 122 &  87 & 1902 &  2.09 \\
79 & 7882 & 131 &  93 & 1894 &  2.24 \\
83 & 7859 & 124 & 101 & 1916 &  2.25 \\
89 & 7878 & 105 &  84 & 1933 &  1.89 \\
97 & 7881 & 106 &  97 & 1916 &  2.03 \\
\hline
\end{tabular}
\end{table}

\begin{table}[htbp!]
\centering
\footnotesize
\setlength{\tabcolsep}{4pt} 
\caption{Agreement and disagreement between the twist model and proxy model (k=7).}
\label{tab:disagreement_updated}
\begin{tabular}{cccccc}
\hline
\textbf{Prime} & \textbf{\makecell{Both \\ Correct}} & \textbf{\makecell{Twist Correct \\ Proxy Wrong}} & \textbf{\makecell{Proxy Correct \\ Twist Wrong}} & \textbf{\makecell{Both \\ Wrong}} & \textbf{\makecell{Disagreement \\ (\%)}} \\
\hline
2  & 6949 & 517 & 525 & 2009 & 10.42 \\
3  & 6613 & 466 & 487 & 2434 &  9.53 \\
5  & 7447 & 307 & 289 & 1957 &  5.96 \\
7  & 7682 & 229 & 242 & 1847 &  4.71 \\
11 & 7837 & 198 & 195 & 1770 &  3.93 \\
13 & 7877 & 196 & 178 & 1749 &  3.74 \\
17 & 7834 & 216 & 178 & 1772 &  3.94 \\
19 & 7864 & 157 & 165 & 1814 &  3.22 \\
23 & 7883 & 166 & 191 & 1760 &  3.57 \\
29 & 7884 & 160 & 143 & 1813 &  3.03 \\
31 & 7870 & 137 & 152 & 1841 &  2.89 \\
37 & 7883 & 134 & 130 & 1853 &  2.64 \\
41 & 7863 & 146 & 163 & 1828 &  3.09 \\
43 & 7872 & 120 & 127 & 1881 &  2.47 \\
47 & 7863 & 140 & 142 & 1855 &  2.82 \\
53 & 7894 & 144 &  98 & 1864 &  2.42 \\
59 & 7877 & 129 & 125 & 1869 &  2.54 \\
61 & 7897 & 109 &  97 & 1897 &  2.06 \\
67 & 7863 & 111 & 114 & 1912 &  2.25 \\
71 & 7876 & 123 & 101 & 1900 &  2.24 \\
73 & 7887 & 124 &  94 & 1895 &  2.18 \\
79 & 7880 & 133 &  83 & 1904 &  2.16 \\
83 & 7848 & 135 &  97 & 1920 &  2.32 \\
89 & 7869 & 114 &  82 & 1935 &  1.96 \\
97 & 7868 & 119 & 105 & 1908 &  2.24 \\
\hline
\end{tabular}
\end{table}

\begin{algorithm}[htbp!]
\caption{Sign-based twist class matching classifier for predicting $a_p(E)$}
\label{alg:sign_matching_j}
\begin{algorithmic}[1]

\Require
Training set $X$ and test set $T$ of elliptic curves with good reduction at a prime $p<100$. The input in a row $i$ is a sequence of the traces of Frobenius $X_i=(a_q(E))_{q<100, q \ne p}$ and twist hash $h(E)$. The output is $a_p(E_i)$.
\Ensure
Predictions $a_p(E)$ on the test set.

\State \textbf{Step 0 (Precomputation).}
\State Compute the probability distribution $\mathcal{P}(p)$ of $a_p(E)$ on the training set $X$.
\State \textbf{Step 1 (Sign-Pattern Indexing).}
\For{all $E \in X$}
\State Compute the pair $\big(\operatorname{sgn}(E, p),\, \operatorname{sgn}(a_p(E))\big)$, where 
$$
\operatorname{sgn}(E, p) = \big(\operatorname{sgn}(a_2(E)),\, \operatorname{sgn}(a_3(E)),\, \dots,\, \widehat{\operatorname{sgn}(a_p(E))},\, \dots,\, \operatorname{sgn}(a_{97}(E))\big),
$$ \indent where we exclude $\operatorname{sgn}(a_p(E))$ from the tuple $\operatorname{sgn}(E)$.
\EndFor
\State Initialize an empty list $L$ for lookup.
\For{each distinct twist hash $h$ in $X$, and for each pair of curves $(E_1, E_2)$ with $h(E_1)=h(E_2)=h$,}
\State Compute the pair $\big(\operatorname{sgn}(E_1, E_2),\, \operatorname{sgn}_p(E_1, E_2)\big)$ where
\[
\operatorname{sgn}(E_1, E_2) = \operatorname{sgn}(E_1, p) \ast \operatorname{sgn}(E_2, p), \]
\indent $\ast$ denoting pointwise multiplication, and 
\[\operatorname{sgn}_p(E_1, E_2) = \operatorname{sgn}(a_p(E_1)) \operatorname{sgn}(a_p(E_2)).
\]

\State Store the pair of $\big(\operatorname{sgn}(E_1, E_2),\, \operatorname{sgn}_p(E_1, E_2)\big)$ in the list $L$.
\EndFor

\State \textbf{Step 2 (Matching Twist Hash and Majority Voting).}
\For{each test curve $E_T \in T$}
\If{$h(E_T)$ is not found in $X$}
\State Predict $a_p(E_T)$ based on the probability distribution $\mathcal{P}(p)$.  
\Else
\State Initialize a map $\freq: \Z \rightarrow \Z_{\ge 0}$ with $\freq[n]=0$ for all $n \in \Z$.
\State Compute $\operatorname{sgn}(E_T, p)= \big(\operatorname{sgn}(a_2(E_T)),\, \operatorname{sgn}(a_3(E_T)),\, \dots,\, \widehat{\operatorname{sgn}(a_p(E_T))},\, \dots,\, \operatorname{sgn}(a_{97}(E_T))\big).$
\For{each curve $E_X \in X$ with $h(E_X)=h(E_T)$}
\State Compute $\big(\operatorname{sgn}(E_T, E_X),\, a_p(E_X)\big)$.
\If{$\operatorname{sgn}(E_T, E_X)$ occurs in the list $L$}
\State Find all pairs $(E_i, E_j)$ such that
$\operatorname{sgn}(E_T, E_X) = \operatorname{sgn}(E_i, E_j).$
\For{each pair $(E_i, E_j)$,}
\State Compute $s=\sgn_p(E_i, E_j)$.
\State Update $\freq[s \cdot a_p(E_X)] \leftarrow \freq[s \cdot a_p(E_X)] + 1$.
\EndFor
\EndIf
\EndFor
\If{$\freq$ is the constant zero map,} 
\State Predict $a_p(E_T)$ by sampling from the empirical distribution of $a_p(E_X)$ over $\{E_X \in X \mid h(E_X) = h(E_T)\}$.
\Else
\State Predict $a_p(E_T)=n$ where $\freq[n]$ is maximal. If there are multiple such $n$, predict $a_p(E_T)$ by sampling from the tied candidates according to $\mathcal{P}(p)$.
\EndIf
\EndIf
\EndFor
\end{algorithmic}
\end{algorithm}

\begin{algorithm}[htbp!]
\caption{Sign-based matching classifier using proxy keys for predicting $a_p(E)$}
\label{alg:sign_matching_proxy}
\small
\begin{algorithmic}[1]

\Require
Training set $X$ and test set $T$ of elliptic curves with good reduction at a prime $p<100$. The input in a row $i$ is a sequence of the traces of Frobenius $X_i=(a_q(E))_{q<100, q \ne p}$. The output is $a_p(E_i)$.
\Ensure
Predictions $a_p(E)$ on the test set.

\State \textbf{(Compute proxy keys for all curves).}
\For{all curves $E \in X \cup T$}
\State Compute the proxy key $k(E)$ using the absolute values of the traces at the $k$ largest primes:
$$
k(E) = \big(|a_{p_1}(E)|, \dots, |a_{p_k}(E)|\big)
$$
\EndFor

\State \textbf{Step 0 (Precomputation).}
\State Compute the probability distribution $\mathcal{P}(p)$ of $a_p(E)$ on the training set $X$.
\State \textbf{Step 1 (Sign-Pattern Indexing).}
\For{all $E \in X$}
\State Compute the pair $\big(\operatorname{sgn}(E, p),\, \operatorname{sgn}(a_p(E))\big)$, where 
$$
\operatorname{sgn}(E, p) = \big(\operatorname{sgn}(a_2(E)),\, \operatorname{sgn}(a_3(E)),\, \dots,\, \widehat{\operatorname{sgn}(a_p(E))},\, \dots,\, \operatorname{sgn}(a_{97}(E))\big),
$$ \indent where we exclude $\operatorname{sgn}(a_p(E))$ from the tuple $\operatorname{sgn}(E)$.
\EndFor
\State Initialize an empty list $L$ for lookup.
\For{each distinct proxy key $k$ in $X$, and for each pair of curves $(E_1, E_2)$ with $k(E_1)=k(E_2)=k$,}
\State Compute the pair $\big(\operatorname{sgn}(E_1, E_2),\, \operatorname{sgn}_p(E_1, E_2)\big)$ where
\[
\operatorname{sgn}(E_1, E_2) = \operatorname{sgn}(E_1, p) \ast \operatorname{sgn}(E_2, p), \]
\indent $\ast$ denoting pointwise multiplication, and 
\[\operatorname{sgn}_p(E_1, E_2) = \operatorname{sgn}(a_p(E_1)) \operatorname{sgn}(a_p(E_2)).
\]

\State Store the pair of $\big(\operatorname{sgn}(E_1, E_2),\, \operatorname{sgn}_p(E_1, E_2)\big)$ in the list $L$.
\EndFor

\State \textbf{Step 2 (Matching Proxy key and Majority Voting).}
\For{each test curve $E_T \in T$}
\If{$k(E_T)$ is not found in $X$}
\State Predict $a_p(E_T)$ based on the probability distribution $\mathcal{P}(p)$.
\Else
\State Initialize a map $\freq: \Z \rightarrow \Z_{\ge 0}$ with $\freq[n]=0$ for all $n \in \Z$.
\State Compute $\operatorname{sgn}(E_T, p)= \big(\operatorname{sgn}(a_2(E_T)),\, \operatorname{sgn}(a_3(E_T)),\, \dots,\, \widehat{\operatorname{sgn}(a_p(E_T))},\, \dots,\, \operatorname{sgn}(a_{97}(E_T))\big).$
\For{each curve $E_X \in X$ with $k(E_X)=k(E_T)$}
\State Compute $\big(\operatorname{sgn}(E_T, E_X),\, a_p(E_X)\big)$.
\If{$\operatorname{sgn}(E_T, E_X)$ occurs in the list $L$}
\State Find all pairs $(E_i, E_j)$ such that
$\operatorname{sgn}(E_T, E_X) = \operatorname{sgn}(E_i, E_j).$
\For{each pair $(E_i, E_j)$,}
\State Compute $s=\sgn_p(E_i, E_j)$.
\State Update $\freq[s \cdot a_p(E_X)] \leftarrow \freq[s \cdot a_p(E_X)] + 1$.
\EndFor
\EndIf
\EndFor
\If{$\freq$ is the constant zero map,} 
\State Predict $a_p(E_T)$ by sampling from the empirical distribution of $a_p(E_X)$ over $\{E_X \in X \mid k(E_X) = k(E_T)\}$.
\Else
\State Predict $a_p(E_T)=n$ where $\freq[n]$ is maximal. If there are multiple such $n$, predict $a_p(E_T)$ by sampling from the tied candidates according to $\mathcal{P}(p)$.
\EndIf
\EndIf
\EndFor
\end{algorithmic}
\end{algorithm}

\end{document}